\documentclass[a4paper,12pt,fleqn]{article}

\usepackage{amsmath,amssymb} 
\usepackage{dsfont} 
\usepackage{enumerate}
\usepackage{a4wide}
\usepackage{WCERS-abbreviations}
 
\DeclareSymbolFont{rsfs}{U}{rsfs}{m}{n}       
\DeclareSymbolFontAlphabet{\mathrsfs}{rsfs}  

\title{Weak convergence of the empirical process and the rescaled empirical distribution function in the Skorokhod
     product space}
\author{Dietmar Ferger \and Daniel Vogel\footnote{Supported by the German Research Foundation (Collaborative Research Center 475 at the Dortmund Institute of Technology).}}

\begin{document}
\maketitle

\begin{abstract}
We prove the asymptotic independence of the empirical process $\al_n = \sqrt{n}( \mathds{F}_n - F)$ and the rescaled empirical distribution function $\ba_n = n (\dsF_n(\tau+\frac{\cdot}{n})-\dsF_n(\tau))$, where $F$ is an arbitrary cdf, differentiable at some point $\tau$, and $\dsF_n$ the corresponding empricial cdf. This seems rather counterintuitive, since, for {\it every} $n \in \N$, there is a deterministic correspondence between $\al_n$ and $\ba_n$.

Precisely, we show that the pair $(\al_n,\ba_n)$ converges in law to a limit having independent components, namely a time-transformed Brownian bridge and a two-sided Poisson process. Since these processes have jumps, in particular if $F$ itself has jumps, the Skorokhod product space $D(\R) \times D(\R)$ is the adequate choice for modeling this convergence in. We develop a short convergence theory for $D(\R) \times D(\R)$ by establishing the classical principle, devised by Yu.\ V.\ Prokhorov, that finite-dimensional convergence and tightness imply weak convergence. Several tightness criteria are given. Finally, the convergence of the pair $(\al_n,\ba_n)$ implies convergence of each of its components, thus, in passing, we provide a thorough proof of these known convergence results in a very general setting. In fact, the condition on $F$ to be differentiable in at least one point is only required for $\ba_n$ to converge and can be further weakened.
\end{abstract}

\rm


\section{Introduction} 

This paper brings together two important convergence results in empirical process theory. The first one is the convergence in law of the uniform empirical process (u.e.p.) $\sqrt{n} ( \mathds{G}_n(t) - t )$, $ t \in [0,1]$, to the Brownian bridge. Here $\mathds{G}_n$ denotes the uniform empirical distribution function (u.e.d.f). This result is originally due to M.\ D.\ Donsker \cite{Donsker1952}, who carried out an idea by J.\ L.\ Doob \cite{Doob1949}. 
The work was motivated by the pioneer papers of A.\ N.\ Kolmogorov \cite{Kolmogorov1933} and N.\ V.\ Smirnov \cite{Smirnov1944} about the limit distribution
of $\sup_{t \in [0,1]}|\sqrt{n} ( \mathds{G}_n(t) - t )|$ and $\sup_{t \in [0,1]}\sqrt{n} ( \mathds{G}_n(t) - t )$, respectively.

The other one is the convergence of the {\it rescaled uniform empirical distribution function} (r.u.e.d.f.) $n \mathds{G}_n(\frac{t}{n})$, $t \ge 0,$  to the Poisson process having intensity $1$. Although being nowadays a standard exercise in empirical process theory, the origin of this result has remained, up to this day, unknown to us. It appears in different levels of generality e.g.\ in \cite{Kabanov1980}, \cite{Al-Hussaini1984}  or \cite{Csorgo1988}. 

The Brownian bridge, closely linked to the Brownian motion, and the Poisson process are two fundamental stochastic processes, the relevance of which goes far beyond being limit processes in asymptotic statistics. The empirical distribution function and derived processes (such as the empirical process) are an important field of study in mathematical statistics. See for example \cite{Shorack1986} or \cite{Vaart1996} for a profound treatment of up-to-date empirical process theory with particular focus on statistical applications.

\bigskip
The aim of this paper is to prove the asymptotic independence of the u.e.p.\ and the r.u.e.d.f.,\ but we are going to do so in a general setting. Instead of being uniformly distributed, we let the underlying sequence of i.i.d.\ random variables $\{ X_n \}$ be sampled from an arbitrary distribution function $F$. Then we look at the following generalizations of the u.e.p.\ and r.u.e.d.f., respectively,
\[
        \al_n^F(t) = \sqrt{n} \big( \mathds{F}_n(t) - F(t)\big), \qquad t \in \R,
\] \[
    \ba_n^{F,\tau}(t)   =  \begin{cases}
                                        n \big[ \dsF_n(\tau+\frac{t}{n})-\dsF_n(\tau)\big],  & \mbox{ if } t \ge 0, \\
                                        n \big[ \dsF_n(\tau+\frac{t}{n})-\dsF_n(\tau-)\big], & \mbox{ if } t < 0,
                                    \end{cases} 
\]
where $\tau$ is an arbitrary real constant and
\[
\mathds{F}_n(t) = \frac{1}{n} \sum\limits_{k = 1}^{n} \Ind{ X_k \le t}, \qquad t \in \R
\]
is the empirical cdf corresponding to $F$. The processes $\AL$ and $\BA$ both converge in law - to limits, say $B_1$ and $N_0$, respectively, that will be properly specified in section \ref{mainresult}. We are going to show that also
\be \label{mainresult-kurz}
    (\AL, \BA) \cid (B_1,N_0),
\ee
where $B_1$ and $N_0$ are stochastically independent.

\bigskip
At this point we would like to spare a few words about the implications of (\ref{mainresult-kurz}). It is quite a remarkable result. The fact that the convergence extends from the individual sequences to the joint sequence is, although not to be taken for granted, hardly surprising. But $B_1$ and $N_0$ in (\ref{mainresult-kurz}) are independent, while $\AL$ and $\BA$ -- since derived from the same sequence $\{X_n\}$ -- are clearly not. Consider a fixed $t \in \R$. One implication of (\ref{mainresult-kurz}) is that $\AL(t)$ and $\BA(t)$ are asymptotically independent. This may seem plausible, since they are deterministic transformations of $\dsF_n(t)$ and $\dsF_n(\frac{t}{n})$, respectively, and it is known that the extreme and middle order statistics are asymptotically independent, cf.\ \cite{Rossberg1967}\footnote{The authors thank E. H\"ausler for pointing out the reference.}. But (\ref{mainresult-kurz}) states even stronger that the whole processes are asymptotically independent -- and that although, for any fixed $n$, $\al_n$ and $\ba_n$ are linked via the strongest form of stochastic dependence there is: knowing one means knowing the other.

\bigskip
When it comes to proving the result, the first question arising is: weak convergence in which measurable space? Since we canonically take Borel-$\sigma$-fields, it comes down to choosing a topological space, which desirably is metrizable and separable. 
The processes involved have discontinuous paths and the whole real line as their time domain. Thus for example the nice, separable metric space $(C[0,1],||\cdot||_{\infty})$, the space of all continuous functions on $\I$, is not an option. But the trajectories of all processes are right-continuous, and the left-hand limits exist in all points, i.e.\ they are  {\it c\`adl\`ag} functions: ``continue \`a droite, limites \`a gauche'' (sometimes also {\it rcll}). The space of all c\`adl\`ag functions on the time domain $T$ is usually denoted by $D(T)$.

An element of $D(\R)$ stays bounded on a compact set, just as a continuous function does. Hence the sup-metric $||\cdot||_{\infty}$ is a possible metric for $D[0,1]$. It induces the {\it topology of uniform convergence} or short, the {\it uniform topology}. However, this metric is unsuitable for $D[0,1]$, due to several reasons. First, $(D[0,1],||\cdot||_{\infty})$ is not separable (see e.g. \cite{Jacod2002}, page 325). Second, and more severe, there are measurability problems. The empirical process is not measurable with respect to the uniform topology. In fact, Donsker's original proof of the weak convergence of the u.e.p. was flawed, because he used this topology.

In 1956, A.\ V.\ Skorokhod \cite{Skorokhod1956} proposed several other topologies on $D[0,1]$, of which the $J_1$-topology has become the most popular. It is coarser than the uniform topology, separable, metrizable and solves the measurability issue. It allows for a workable Arzela-Ascoli-type compactness characterization and it also declares a convergence more natural to functions with jumps. Nowadays, $D[0,1]$ is by default equipped with the $J_1$-topology and simply referred to as the {\it Skorokhod space}. We endow $D(\R)$ with a proper extension of this $J_1$-topology (by the same means one declares a uniform topology for functions on the real line, cf. page \pageref{cc1}) and then treat the convergence statement (\ref{mainresult-kurz}) in the Skorokhod product space $D(\R) \times D(\R)$. The proof then breaks down into two tasks: Derive a weak convergence criterion in the space $D(\R) \times D(\R)$ (Theorem \ref{wc2}) and show that $(\al_n^F,\ba_n^{F,\tau})$ satisfies it (section \ref{proof}).

The standard method of proving weak convergence of stochastic processes is as follows: Prove the weak convergence of the finite-dimensional distributions, and show that the sequence is tight. The key argument here is Prokhorov's theorem \cite{Prokhorov1956}. For example, this method is used to show that the partial sum process converges to the Brownian motion in $(C[0,1],||\cdot||_{\infty})$ (Donsker's theorem \cite{Donsker1951}). The principle transfers with little alteration to $D[0,1]$ and $D(\R)$. We will show that it extends as well to $D(\R) \times D(\R)$. It is, however, only feasible, if the finite-dimensional distributions are known, and there are other approaches as well, see e.g. \cite{Jacod2002}.

\bigskip
The paper is organized as follows: Section \ref{mainresult} states the task in detail, the principal statement of this paper is formulated in Theorem \ref{main result}. The predominant 
rest of the paper is devoted to its proof: Section \ref{WCinD} introduces the space $D(\R)$ and states the classic convergence 
criterion, section \ref{TinD} deals with tightness in $D(\R)$. 
In section \ref{WCinDtimesD} we begin to develop a short weak convergence theory for the product space $D(\R) \times D(\R)$ and 
prove an analoguous convergence criterion. Finally we apply the latter to show Theorem \ref{main result} in section \ref{proof}.
The paper ends with section 7 in which a short description is given how Theorem 2.1 can be used in statistics.

\bigskip
We conclude the introduction with some remarks on the literature. Most of what we use are classical results being covered in a variety of textbooks. Our main reference is Patrick Billingsley's ``Convergence of Probability Measures'' \cite{Bill99}. This book's first edition dates back to 1968 and features a stage-wise development from $C[0,1]$ to $D[0,1]$ to $D[0,\infty)$. A number of newer books, like e.g.\ \cite{Pollard1984}, \cite{Ethier1986}, \cite{Whitt2002} and \cite{Jacod2002} consider right away the space $D[0,\infty)$, or a more general version of it, without paying extra attention to $D[0,1]$, and hence tend to be more profound. In particular, \cite{Jacod2002} gives an exhaustive treatment of weak convergence on $D[0,\infty)$.

Besides \cite{Skorokhod1956} the other two important papers on $D\I$ are Kolmogorov \cite{Kolmogorov1956} and Prokhorov \cite{Prokhorov1956}. The analogue on $D[0,\infty)$ is due to C. Stone \cite{Stone1963}.
Billingsley \cite{Bill99} gives a construction of a complete metric on $D[0,\infty)$. He adopts a suggestion of T. Lindvall \cite{Lindvall1973}, who in turn follows W. Whitt's approach on $C[0,\infty)$, \cite{Whitt1970}. Whitt also suggests another metric on $D[0,\infty)$, \cite{Whitt1971}.

\section{Main result} \label{mainresult}

Let $F$ be an arbitrary distribution function, and $X_1, X_2, ...$ a sequence of i.i.d.\ random variables being distributed according to $F$. The corresponding empirical distribution function (edf) is given by
\[
     \mathds{F}_n(t) = \frac{1}{n} \sum\limits_{k = 1}^{n} \Ind{ X_k \le t}, \qquad t \in \R,\ n \ge 1.
\]
The following family $\{\AL |\ n \in \N \}$ of random functions is called the {\it empirical process}:
\[
    \al_n^F(t) = \sqrt{n} \big( \mathds{F}_n(t) - F(t)\big) = \frac{1}{\sqrt{n}} \sum_{k=1}^n \big( \ind{X_k \le t} - F(t) \big),
    \qquad t \in \R, n \ge 1.
\]
Furthermore, for any real number $\tau$, let
\begin{eqnarray} 
 \ba_n^{F,\tau}(t)  &=&  \begin{cases}
                                        n \big[ \dsF_n(\tau+\frac{t}{n})-\dsF_n(\tau)\big],  & \mbox{ if } t \ge 0, \\
                                        n \big[ \dsF_n(\tau+\frac{t}{n})-\dsF_n(\tau-)\big], & \mbox{ if } t < 0,
                                    \end{cases}    \nonumber \\
                            & =&  \begin{cases}
                                        \sum_{k=1}^n \Varind{(\tau, \tau+\frac{t}{n}]} (X_k), & \mbox{ if } t \ge 0, \\
                                        \sum_{k=1}^n -\Varind{(\tau+\frac{t}{n},\tau)} (X_k), & \mbox{ if } t < 0.
                                    \end{cases} \label{alpha^F,tau}
\end{eqnarray}
We want to call the family $\{\ba_n^{F,\tau} |\ n \in \N \}$ the {\it rescaled empirical distribution function}. Whenever it is clear or not of interest which $F$ and $\tau$ are meant, we will shortly write $\al_n$ and $\ba_n$. In Theorem \ref{main result} we will make the following basic assumptions on $F$.

\medskip
{\bf Condition C.1} $F$ has both, left- and right-hand side, derivatives in $\tau$. Call the former $\varrho_1$ and the latter $\varrho_2$, i.e.\
\be \label{rho_1}
    \varrho_1 = \lim_{h \nearrow 0}\frac{F(\tau+h)-F(\tau-)}{h} 
\ee
and
\be
    \varrho_2 = \lim_{h \searrow 0}\frac{F(\tau+h)-F(\tau)}{h}. 
\ee
Pay attention to the $\tau-$ in line (\ref{rho_1}). This definition of left-hand side derivative does not require $F$ to be continuous in $\tau$. 

\medskip
The next step is to specify the limit processes of $\{ \al_n \}$ and $\{ \ba_n \}$. Let $B_0 = \{B_0(t) |\, t \in \I\}$ be a Brownian bridge and
\[
    B_1 = B_1^F = B_0 \circ F,\qquad \mbox{i.e. } \quad B_1^F(t) = B_0(F(t)), \qquad t \in \R.
\]
$B_1^F$ is a Gaussian process with expectation zero and covariance function $\cov(s,t) = F(s)(1-F(t))$ for $s \le t$.
Furthermore, let $N_1$, $N_2$ be two Poisson processes with the following properties:
\begin{itemize}
\item
    $N_1$ and $N_2$ are independent of $B_1$, and of each other.
\item
    $N_i$ has rate $\varrho_i$, $i = 1,2$.
\item
    $N_2$ has, as usual, right-continuous trajectories while those of $N_1$ are left-continuous, i.e.\ the value at a jump point is always set to the left-hand limit. Note that this leaves the finite-dimensional distributions unchanged.
\end{itemize}
Then define $N_0 = N_0^{\varrho_1,\varrho_2}$ by
\[
 N_0(t) = \begin{cases}
                            -N_1(-t), &  \qquad t < 0,  \\
                             N_2(t),  &  \qquad t \ge 0.
                        \end{cases}
\]
All four stochastic processes we have introduced so far, $\al_n$, $\ba_n$, $B_1$ and $N_0$, have trajectories in the Skorokhod space $D = D(\R)$, that is the space of all \it c\`adl\`ag \rm functions on $\R$, equipped with the Skorokhod topology ($J_1$-topology). The Skorokhod space is properly introduced in section \ref{WCinD}. Whenever we write $D$, the topological space is meant. This also applies to $D \times D$ (product topology). The Skorokhod space $D$ is a Polish space. Its Borel-$\sigma$-field shall be denoted by $\D$. 

It is known and in the case of the uniform(0,1) distribution considered to be folklore that
\begin{enumerate}[(A)]
\item $\AL \cid B_1^F$ in $D$, and
\item \label{II} if C.1 is satisfied, $\BA \cid \Nrho$ in $D$.
\end{enumerate}

{\bf Remark on (\ref{II}).} 
The definition (\ref{alpha^F,tau}) of $\{\ba_n\}$ 
resembles a
sequence of difference quotients. 
It is therefore not suprising that the derivative of $F$ at $\tau$ appears as parameter in the limit process of $\{\beta_n\}$. But $F$ does not have to be differentiable at $\tau$, it may have a sharp bend at $\tau$ (i.e.\ left- and right-hand derivative differ) or even a jump (left- and right-hand limit differ). The behavior of the process $\ba_n$ on the positive half-axis is determined by the behavior of $F$ in the right-hand vicinity of $\tau$, likewise for the negative half-axis. Thus the assumption of differentiability can be weakend to Condition C.1.

In the proof of Proposition \ref{al,ba fidis} it becomes clear that, instead of C.1, we only require the (formally) weaker condition that all of the limits 
\[ \textstyle
    \lim\limits_{n \to \infty} \frac{n}{t} \big[ F(\tau+\frac{t}{n})-F(\tau-)\big] \mbox{ for } t < 0 \quad \mbox{and} \quad
    \lim\limits_{n \to \infty} \frac{n}{t} \big[ F(\tau+\frac{t}{n})-F(\tau)\big] \mbox{ for } t > 0
\]
exist in order to have convergence of the process $\{ \ba_n \}$. However, it can be shown that for monotone functions $F$ these limits coincide for either all $t > 0$ or all $t < 0$ and, even more, that the left- and right-hand derivatives exist at $\tau$. Finally it should be noted that C.1 is a very weak condition  -- it does, for instance, not imply continuity in an upper or lower neighborhood of $\tau$. 

\medskip
The following result is new.

\bt \label{main result} 
Under  Condition C.1 \ $(\al_n^F,\ba_n^{F,\tau})$ converges in distribution to $(B_1^F,N_0^{\varrho_1,\varrho_2})$ in $D \times D$.
\et

\bf Remarks. \rm
\begin{enumerate}[(I)]
\item 
Keep in mind that we have defined $N_0$ to be independent of $B_1$, which specifies the distribution of $(B_1,N_0)$. The remarkable feature of Theorem \ref{main result} is not the convergence itself, but rather the fact, that the ``highly  dependent'' $\al_n$ and $\ba_n$ (knowing one means knowing the other) converge to independent limits.
\item
Of course, (\ref{II}) follows from Theorem \ref{main result}. Note that -- apart from the regularity condition C.1 -- $F$ is completely arbitrary. The result does not seem to be contained as such in the literature. It should, however, be compared to Theorem 3.1.\ in \cite{Csorgo1988}. The authors there consider processes of the type $n \dsF_n(a_n t + b_n)$, where $\{a_n\}$ and $\{b_n\}$ are sequences of real numbers such that the adjusted first order statistic $(X_{1:n}- b_n)/a_n$ converges to a non-degenerate limit. Such an extreme-value process may coincide with $\BA$ if $\tau$ is the left endpoint of the support of $F$ and $\{b_n\}$ is constant equal to $\tau$. In this situation, our condition C.1 with $\varrho_2 > 0$ implies the assumption on $F$ in \cite{Csorgo1988}: $F$ lies in the domain of attraction of the cdf $L_{2,1}(x) = (1-e^{-x})\varind{(0,\infty)}(x)$ (Weibull distribution with shape parameter 1), which is an extreme-value distribution of type 2, cf.\ e.g.\ \cite{Galambos1978}, pp.\ 58, 76.
\item 
Convergence in law is canonically defined on a Borel-$\sigma$-field, the underlying topological space being here $D \times D$. On the other hand, $(\al_n, \ba_n)$ and $(B_1, N_0)$ are pairs of random variables and hence defined on the product measure space, i.e.\ the $\sigma$-field $\D \otimes \D$. Fortunately, $\D \otimes \D$ coincides with the Borel-$\sigma$-field on $D \times D$ (cf.\ Lemma \ref{sigma-fields}).
\end{enumerate}

The proof of Theorem \ref{main result} is subject of section \ref{proof}. 

\section{Weak Convergence in \boldmath$D$} \label{WCinD}

Preliminary note: Most textbooks and articles consider the space $D[0,\infty)$ instead of $D(\R)$. Both spaces are qualitatively equal, all results for $D[0,\infty)$ hold with little notational change (which demands its due amount of care) also for $D(\R)$. Define
\[
    D = D(\R) = \{ x:\R \to \R\ |\ x(t-),\ x(t+) \mbox{ exist}, \ x(t) = x(t+) \ \ \forall \ t \in \R \}.
\]
Elements of $D$ have at most countably many discontinuity points and are bounded on compact sets. We declare a topology on $D$ \rm by the following characterization of convergence. Let $\Lambda$ denote the class of all strictly increasing, continuous, surjective mappings $\lambda$ from $\R$ \emph{onto} itself. A sequence $\{ x_n \} \subset D$ converges to $x \in D$ if and only if a sequence $\{\lambda_n\} \subset \Lambda$ exists such that
\be \label{cc1} 
        \begin{cases} 
            \lambda_n(t) \to t              & \qquad \mbox{ uniformly in } t \in \R, \\
            x_n(\lambda_n(t)) \to x(t)  & \qquad \mbox{ uniformly in } t \in [-m,m] \mbox{ for all } m \in \N.  
        \end{cases}
\ee

This is a $D(\R)$-version of the {\it $J_1$-topology}, originated by  A.\ V.\ Skorokhod \cite{Skorokhod1956}. This is the only topology we consider on $D$ and subsequently refer to it as the \it Skorokhod topology\rm. For details on different topologies on $D$ see for example \cite{Revuz1999}. Compare the above characterization to 1.14, page 328, in \cite{Jacod2002}. Note that, unlike in $D[0,\infty)$, the point $0$ must not play a special role in $D(\R)$.

Since the identity is an element of $\Lambda$, uniform convergence on compact sets implies Skorokhod convergence. In fact, the Skorokhod topology is strictly coarser than the topology of locally uniform convergence (\it uniform topology\rm ), i.e.\ there are fewer open sets and more convergent sequences. For instance, 
$\{ \varind{[\frac{1}{n},\infty)} |\ n \in \N \}$ convergences in the Skorokhod topology, but not in the uniform topology. As mentioned before, by writing $D$ we always mean the topological space. Let $\D$ be its Borel-$\sigma$-field. The topological space $D$ is separable (whereas the set $D$ endowed with the uniform topology is not separable), completely metrizable and in this sense a Polish space. See e.g.\ \cite{Bill99} for a complete metric.

\medskip
Let $\pi_t$ denote the projection $\pi_t: D \to \R: x \mapsto x(t)$. For any probability measure $P$ on $(D,\D)$ let $T_P$ be the set of all points $t \in \R$ for which $\pi_t$ is $P$-almost surely continuous. We call $P_X$ the distribution of any random variable $X$ in $D$ and write $T_X$ for $T_{P_X}$. 

\bl \label{camc} 
The complement of\, $T_X$ in $\R$ is at most countable. 
\el  

See e.g.\ \cite{Bill99}, page 174. It follows, that $T_X$ is dense.

\bd \label{cofidis} 
In $D$ we say, \it the finite-dimensional distributions (fidis) of $X_n$ converge to those of $X$,\rm\ and write $X_n \fdc X,$ if 
\be \label{fidis}
    \big( X_n(t_1), ..., X_n(t_k) \big) \cid \big( X(t_1), ..., X(t_k) \big)     
\ee     
for all $k \in \N$ and $t_1 < ... < t_k \in T_X$.
\ed

\bf Remarks. \rm
\begin{enumerate}[(I)]
\item 
We restrict $t_1,...,t_k$ to lie in  $T_X$, because $X_n \cid X$ in $D$ does not necessarily imply $\pi_t(X_n) \cid \pi_t(X)$, cf.\ e.g.\ \cite{Jacod2002}, page 349, 3.14. This is due to the fact that $\pi_t$ (function from $D$ to $\R$) is continuous at a point $x$ only if $x$ (function from $\R$ to $\R$) is continuous at $t$, cf. \cite{Bill99}, page 134, Theorem 12.5 (i). Think, for instance, of $\varind{[\frac{1}{n},\infty)} \ra \varind{[0,\infty)}$, but $\pi_t(\varind{[\frac{1}{n},\infty)}) \nrightarrow \pi_t(\varind{[0,\infty)})$.
\item
\ref{cofidis} is equivalent to: there exists a \it dense \rm subset $S$ of $\R$ such that (\ref{fidis}) holds for all finite subsets $\{ t_1,...,t_k \}$ of $S$, see \cite{Jacod2002}, page 350, 3.19. In this sense $\fdc$ does not depend on its right-hand side.
\end{enumerate}

Now here is a characterizations of weak convergence in $D$. It is phrased in terms of random variables and convergence in law - which is equivalent to the weak convergence of the respective distributions. 

\bp \label{wc1}
Let $\{ X_n \}$ be a sequence of random variables in $(D,\D)$ with the following two properties:
\begin{enumerate}[(1)]
    \item $\{ X_n \}$ is tight.
    \item $X_n \fdc X$.
\end{enumerate}
Then $X_n \cid X$.
\ep

Short, convergence of the fidis (in the sense of \ref{cofidis}) and tightness together imply convergence in law. These two conditions are sufficient and necessary, cf.\ e.g.\ \cite{Jacod2002}, page 350 or \cite{Bill99}, page 139.

\section{Tightness in \boldmath$D$} \label{TinD}

In order to make use of Proposition \ref{wc1} we need a handy tightness criterion. Recall tightness of a sequence: A family $\mathrsfs{P}$ of probability measures on the Borel-$\sigma$-field of a metric space is {\it tight}, if for every $\varepsilon > 0$ there exists a compact set $K$ such that $P(K) > 1 - \varepsilon$ for every $P \in  \mathrsfs{P}$. A family of random variables is {\it tight} if the family of their respective distributions is tight.
Prokhorov's theorem tells us that in complete metric spaces, like $D$, tightness is equivalent to relative compactness. A family of $\mathrsfs{P}$ of probability measures is \it relatively compact\rm, if every sequence in $\mathrsfs{P}$ contains a convergent subsequence. The limit needs not to lie in $\mathrsfs{P}$.

We present three criteria which allow to confirm that a given sequence of random variables in $D$ is tight. 
The first is, in fact, a characterization of tightness.

\subsection{A tightness characterization in \boldmath$D$}

First we need to introduce some notation. We will be dealing with intervals of the type $[-m,m]$, where $m \in \N$. 

For an arbitrary function $x: \R \ra \R$ and an arbitrary set $T \subset \R$ we define
\be \label{w_x}
    w(x,T) = \sup_{s,t \in T} |x(s)-x(t)|.
\ee
We want to call any finite set $\sigma = \{ s_0, ..., s_{k} \} \subset \R$ satisfying $-m = s_0 < s_1 < ... < s_{k} = m$ a {\it grid on $[-m,m]$}. If 
\[
    s_i - s_{i-1} > \de     \qquad \mbox{for all }  i = 2, ..., k-1,
\]
i.e.\ all intervalls except those at the left and right end are wider than $\de$, we want to call the grid {\it $\de$-sparse}. Let $\Ss(m,\de)$ be the set of all $\de$-sparse grids on $[-m,m]$ and define the following modulus:
\be \label{w.hat}
    \hat{w}_m(x,\de) =  \inf_{\Ss(m,\de)} \max_{1 \le i \le k} w(x,[s_{i-1}, s_i)).
\ee

\bt \label{tc1}
A sequence of random variables $\{ X_n \}$ in $(D,\D)$ is tight if and only if the following two conditions hold.
\begin{enumerate}[(1)]
\item \label{tc1b1}
For all $t$ in a dense subset $T_0$ of $\R$,
\[ 
    \lim_{a \rightarrow \infty} \limsup_n \Pb(|X_n(t)| \ge a) = 0,
\]
\item  \label{tc1b2}
and for every $m \in \N$ and $\varepsilon > 0$,
\[
    \lim_{\de \rightarrow 0} \limsup_n \Pb(\hat{w}_m(X_n,\de) \ge \varepsilon) = 0.
\]
\end{enumerate}
\et

\bop cf.\ \cite{Bill99}, Theorem 16.8 in combination with the subsequent corollary.
\eop

\subsection{A moment-type tightness criterion}

\bp \label{tc2}
Let $X$ and $X_n$, $n \in \N$, be random variables in $\DD$. Suppose that
\begin{enumerate}[(1)]
\item \label{tc2b1}
    $X_n \fdc X,$ and
 \item \label{tc2b2}
    there exists a non-decreasing, continuous function $H:\R \ra \R$ and real numbers $a > 1$ and $b \ge 0$, such that
    \[ 
        \E \Big( |X_n(s) - X_n(r)|^b |X_n(t) - X_n(s)|^b \Big) \le \big(H(t) - H(r) \big)^a.
    \]
    holds for all $r < s < t$, and $n \ge 1$.
\end{enumerate}
Then $\{X_n \}$ is tight. 
\ep

\bop The $D[0,1]$ version of this proposition is Theorem 13.5 on page 142 in \cite{Bill99}. The proof is also worked out in detail for the $D(\R)$-case in \cite{Vogel2005}.

It needs to be shown that \ref{tc2} (\ref{tc2b1}), (\ref{tc2b2}) imply \ref{tc1} (\ref{tc1b1}), (\ref{tc1b2}). In fact, \ref{tc1} (\ref{tc1b1}) follows from \ref{tc2} (\ref{tc2b1}) already, and \ref{tc1} (\ref{tc1b2}) follows from \ref{tc2} (\ref{tc2b2}).
And of course, by Proposition \ref{wc1}, under the assumptions of Proposition \ref{tc2} we have $X_n \cid X$ in $\DD$.
\eop

\subsection{A point-process tightness criterion}

Let $\T$ be the set of all non-decreasing series $\{ t_z |\ z \in \Z \}$ that meet the restrictions $t_z \in [-\infty,\infty]$ for all $z \in \Z$, $t_0 \le 0 < t_1$, $t_z \to \pm \infty$ as $z \to \pm \infty$ and  $\{t_z\}$ is strictly increasing where it is not $\pm \infty$. Then define the following two classes of funtions,
\begin{eqnarray*} 
    \Vp &=& \Big\{ x : \R \to \R\, \Big|\ x = c 
                                                                    - \sum_{z = -\infty}^0 \varind{(-\infty,t_z)}  
                                                                    + \sum_{z =1}^{\infty} \varind{(t_z,\infty)}, \quad c \in \Z, \{t_z\} \in \T \Big\},\\
    \V^+ &=&  \Big\{ x : \R \to \R\, \Big|\ x = c 
                                                                    - \sum_{z = -\infty}^0 c_z \varind{(-\infty,t_z)}  
                                                                    + \sum_{z =1}^{\infty} c_z \varind{(t_z,\infty)}, \\
  & & \hspace{13.0em} c \in \Z, \{t_z\} \in \T, c_z \in \N \mbox{ for all } z \in \Z \ \Big\}. \qquad
\end{eqnarray*}

Apparently $\Vp \subset \V^+ \subset D$. The set $\V^+$ allows also the following characterization: it contains all elements of $D$ that are non-decreasing and integer-valued. Then, by employing (\ref{cc1}), it is easy to see that the potential limit of any series $\{ x_n \}$ in $\V^+$ has this property, too. Hence $\V^+$ is closed in $D$ and therefore measurable. As for $\Vp,$ note that Skorokhod convergence $x_n \to x$  implies that for all $t \in \R$ there exists a sequence $\{t_n\} \subset \R$ such that $t_n \to t$ and 
\[
    x_n(t_n) - x_n(t_n-) \ \lra \ x(t)- x(t-),
\]
cf.\ \cite{Jacod2002}, page 337, 2.1. Hence, if $\{x_n\} \subset \Vp$, the limit $x$ can only have jumps of size 1 as well: $\Vp$ is closed in $D$. We want to call a random function whose paths lie almost surely in $\Vp$ a \emph{counting process}.

\bp \label{tc3}
Let $X$ and $X_n$, $n \in \N$, be random variables in $\DD$. Suppose that
\begin{enumerate}[(1)]
\item \label{tc3b1}
    $X_n \fdc X$,
 \item 
    $\Pb(X \in \Vp) = 1$ and 
 \item
    $\Pb(X_n \in \V^+) = 1$, $n \in \N.$
\end{enumerate}
Then $\{X_n \}$ is tight. 
\ep

This is generalization of Theorem 3.37, page 354, in \cite{Jacod2002}. Basically, \cite{Jacod2002} consider the space $D[0,\infty)$ and require $X_n$, $n \in \N$, also to be counting processes in the above sense. A proposition of exactly the same type as ours ($X$ has jumps of size 1, $X_n$ has integer-valued jumps) can be found in \cite{Csorgo1988}. For the sake of completeness we present an alternative proof.

\medskip
\bopt{ of Proposition \ref{tc3}}
We apply, of course, Theorem \ref{tc1}. The implication \ref{tc3} (\ref{tc3b1}) $\Longrightarrow$ \ref{tc1} (\ref{tc1b1}) is straightforward, cf.\ e.g.\ proof of  Theorem 13.3 in \cite{Bill99}. We only derive condition \ref{tc1} (\ref{tc1b2}) here.

Let $m \in \N$ and initially also $\de > 0$ be fixed. Then choose a $\de$-sparse grid $\sigma = \{ s_0, ..., s_k \}$ on $[-m,m]$ according to the following additional restrictions:
\[
    s_i - s_{i-1} < 2 \de, \qquad i = 1, ..., k, \qquad \mbox{ and } \qquad s_1,...,s_{k-1} \in T_X.
\]
The latter is always possible, since $T_X$ is dense in $\R$, but $s_0 = -m \in T_X$ or $s_k = m \in T_X$ does not need to hold. Now define the following two sets,
\[
    A = \Big\{ (t_1,...,t_{k-1})   \Big|\  t_{i+1} - t_{i-1} < \frac{3}{2}, \quad i = 2,...,k-2  \Big\} \subset \R^{k-1},
\]
\[
    \tAs =  \Big\{ x \in \V^+   \Big|\  \big(x(s_1),...,x(s_{k-1})\big) \in A  \Big\} \subset \V^+.
\]
With these constructions the proof breaks down into two steps. First we show
\begin{enumerate}[(a)]
\item $\limsup\limits_n \Pb(\hat{w}_m(X_n,\de) \ge \varepsilon)  \le  \Pb( X \notin \tAs)$ for all positive $\varepsilon$ and then
\item $\Pb( X \notin \tAs) \to 0$ as $\de \to 0$. 
\end{enumerate}

Part (a): By construction of $\tAs$ we have $x \in \tAs \ \Rightarrow \ \hat{w}_m(x,\de) = 0$, which implies
\[
    \Pb (X_n \notin \tAs) \, \ge  \, \Pb(\hat{w}_m(X_n,\de) \ge \varepsilon) \qquad \forall \ n \in \N,\, \ve > 0.
\]
Since we have chosen $s_i \in T_X$, $i = 1,...,k-1$,
\[
    \Big( X_n(s_1),...,X_n(s_{k-1}) \Big) \cid  \Big( X(s_1),...,X(s_{k-1}) \Big).      
\]
The set $A$ is open in $\R^{k-1}$, and hence by the Portmanteau theorem
\[
    \Pb (X \notin \tAs)\, \ge \, \limsup_n \Pb (X_n \notin \tAs) 
    \, \ge \, \limsup_n\Pb(\hat{w}_m(X_n,\de) \ge \ve) \quad \forall \ \ve > 0.
\]

Part (b): Define $T_z$, $z \in \Z$, to be the jump times of $X$ (we understand them as random variables in $[-\infty,\infty]$), where we count as follows:
\[
    ... < T_{-2} < T_{-1} < T_0 \le 0 < T_1 < T_2 < ...
\]
Now consider the following events
\[
    B_n = \big\{ T_{-n} < -m,\  m < T_n \big\}, \qquad n \in \N, 
\]
\[
    C_{\de,n} = \big\{ T_i - T_{i-1} > 4 \de, \ i = -n+1,...,0,...,n \big\}, \quad n \in \N, \de > 0. 
\]
It holds
\be \label{h1}
    \forall \ \ve > 0 \ \exists \ n \in \N \ :  \ \Pb(B_n) \ge 1 - \frac{\ve}{2},
\ee
\be \label{h2}
    \forall \ \ve > 0, n \in \N \ \exists \ \de > 0 \ :  \ \Pb(C_{\de,n}) \ge 1 - \frac{\ve}{2}
\ee
and
\be \label{h3}
    \{ X \in \Vp \} \cap B_n \cap C_{\de,n} \ \subset \ \{ X \in \tAs \}.
\ee
Since   $\Pb(X \in \Vp) = 1$, (\ref{h1}),  (\ref{h2}) and  (\ref{h3}) imply together 
\[
    \forall \ \ve > 0 \ \exists \ \de > 0 \ : \ \Pb(X \notin \tAs) \le \ve.
\]
It remains to show (\ref{h1}) and (\ref{h2}). Both follow by the same principle from the fact that $X$ is a counting process, which we will exemplify at (\ref{h1}). Assume the opposite is true:
\[
    \exists \ \ve > 0 \ \forall \ n \in \N \ : \ \Pb(B_n) < 1-\frac{\ve}{2}.
\]
Since $\{ B_n |\ n \in \N \}$ is an increasing series of sets, 
\[
    \Pb(\bigcap_n \mathsf{C} B_n) \ge \frac{\ve}{2},
\]
where $\mathsf{C}$ means set complement. The event $\bigcap_n \mathsf{C} B_n$ reads as: all $T_z$, $z \in \Z$, lie in $[-m,m]$. By definition of the set $\T$ this is a contradiction to $X \in \Vp$.
\eop
\section{Weak Convergence in \boldmath$D \times D$} \label{WCinDtimesD}

In Theorem \ref{wc2} we will give a weak convergence characterization in $D \times D$ of the same type as Proposition \ref{wc1}. The set $D \times D$ is the collection of all pairs $(x,y)$, where $x, y \in D$. It is endowed with the product topology, i.e.\ 
\be \label{cc2} 
    (x_n,y_n) \to (x,y) \quad  \Longleftrightarrow \quad
    \begin{cases} 
            x_n \to x, \\
            y_n \to y.  
        \end{cases}
\ee
 Again, by writing $D \times D$ we refer to the topological space. $D \times D$ is a Polish space. The following is important:
\bl \label{sigma-fields}
The Borel-$\sigma$-field on $D \times D$ coincides with the product $\sigma$-field $\D \otimes \D$.
\el
\bop See e.g.\ \cite{Elstrodt2002}, Theorem 5.10, page 115. Separability is needed.
\eop

\bf Remark. \rm One can identify the pair of functions $(x,y)$ with the function $f_{x,y}: \R \to \R^2: t \mapsto (x(t),y(t))$. If and only if $x, y \in D$, then $f_{x,y}$ is a c\`adl\`ag function from $\R$ to $\R^2$. We want to call the space of such functions $D(\R,\R^2)$. The generalization is straightforward: The convergence characterization reads exactly as (\ref{cc1}), only $x_n(\lambda_n(t))$ and $x(t)$ are $\R^2$-valued. In fact the co-domain can easily be replaced by any Polish space without having to change anything.

If we identify $(x,y) \leftrightarrow f_{x,y}$, the sets $D \times D$ and $D(\R,\R^2)$ are equal, but the Skorokhod topology on $D(\R,\R^2)$ is strictly finer than the product topology on $D \times D$, i.e.\ it has less convergent sequences. Take, for instance, $x_n = \varind{[\frac{1}{n},\infty)}$, $y_n = \varind{[-\frac{1}{n},\infty)}$. However, both topologies induce the same Borel-$\sigma$-field, cf.\ \cite{Revuz1999}. In this paper we are not at all concerned with the space $D(\R,\R^2)$. We deal with pairs of random variables and their convergence in law, which we want to be of the same type as (\ref{cc2}). The product topology has to be our concern.

\bt \label{wc2}
Let $\{ Z_n = (X_n, Y_n) \}$ be a sequence of random variables in $(D \times D, \D \otimes \D)$. If
\begin{enumerate}[(1)]
    \item \label{wc2 b1} the sequences $\{ X_n \}$ and $\{ Y_n \}$ are tight, and
    \item \label{wc2 b2} there is a random variable $Z = (X,Y)$ in $\DDDD$ such that
        \be \nonumber 
            \big( X_n(t_1), ..., X_n(t_k), Y_n(t_1), ...,Y_n(t_k) \big)
             \cid           
            \big( X(t_1), ..., X(t_k), Y(t_1), ...,Y(t_k) \big)
        \ee
        for all $k \in \N$, $t_1, ...,t_k \in T_X \cap T_Y$,
\end{enumerate}
\nopagebreak
then $Z_n \cid Z$.
\et 

The rest of the section is devoted to the proof of Theorem \ref{wc2}. It is more convenient to formulate the proof in terms of probability measures than random variables. Therefore, let $P$, $P_n$, $\p{1}$, $\p{1}_n$, $\p{2}$ and $\p{2}_n$ be the distributions of $Z$, $Z_n$, $X$, $X_n$, $Y$ and $Y_n$, $n \in \N$, respectively. 

One thing to note about the theorem is that in \ref{wc2} (\ref{wc2 b1}) we only require $\{ X_n \}$ and $\{ Y_n \}$ individually to be tight. This of course implies tightness of the joint sequence:

\bl \label{Z tight}
If $\{ X_n \}$ and $\{ Y_n \}$ are tight sequences of random variables in $(D,\D)$, then $\{ Z_n = (X_n,Y_n) \}$ is tight in $(D \times D, \D \otimes \D)$. 
\el
\bop
(a corollary of Tikhonov's theorem) For any $\varepsilon > 0$ we find compact sets $K_1, K_2 \subset D$ such that $P^{(i)}_n (K_i) > 1 - \frac{\varepsilon}{2}$\ for all $n \in \N$, $i = 1,2$. Then $P_n(K_1 \times K_2) > 1 - \varepsilon$ for all $n \in \N$, and $K_1 \times K_2$ is compact in $D \times D$ by Tikhonov's theorem.
\eop

We now introduce projections. Let $T = \{t_1,...,t_k\}$ and $S = \{s_1,...,s_l\}$, where $t_1 < ... < t_k$ and $s_1 < ... < s_l$. Define
\[
    \pi_T : D \longrightarrow \R^k : x \mapsto (x(t_1), ..., x(t_k)) = (\pi_{t_1}(x), ...,\pi_{t_k}(x))
\]
and
\[
    \piST : D \times D \ra \R^{l+k}: (x,y) \mapsto (\pi_S(x), \pi_T(y)) = \big( x(s_1), ...,x(s_l),y(t_1), ...,y(t_k) \big).
\]
Then \ref{wc2} (\ref{wc2 b2}) can be written as
\[
    P_n \circ \piTT^{-1} \cid P \circ \piTT^{-1} \qquad \forall\ T \subset T_X \cap T_Y, |T|\ < \infty. 
\]

\bl \label{Paecontinuous}
If \ $T \subset T_X \cap T_Y$, $T$ finite, then $\piTT$ is $P$-a.e.\ continuous.
\el 
\bop
Let $A$ be the discontinuity set of $\pi_T$, i.e.\ the set of all points $x \in D$ in which $\pi_T$ is {\it not} continuous. The function $\piTT$ is continuous at a point $(x,y) \in D \times D$ if and only if $\pi_T$ is continuous at $x$ and $y$. Hence the discontinuity set of $\piTT$ is $(A \times D) \cup (D \times A)$. Due to our assumption $T \subset T_X \cap T_Y$ we have
\[
    P\big( (A \times D) \cup (D \times A)\big) \le P(A \times D) + P(D \times A) = \p{1}(A) + \p{2}(A) = 0.
\]
\vskip-1.5em
\eop

The proof of \ref{wc2} requires furthermore a few measure theoretical concepts. 

\bd \label{sepcla}
Let $(\Omega, \A)$ be a measurable space. Any subclass $\Ss$ of $\A$ that satisfies 
\[
    \mu|_{\Ss} = \nu|_{\Ss}\quad \Rightarrow \quad \mu = \nu
\] 
for any two probability measures $\mu$ and $\nu$ on $\A$ we want to call a \it separating class \rm for $\A$.
\ed
If $\mu$ and $\nu$ differ, then $\Ss$ already suffices to \it separate \rm them. Recall that, if a system of sets $\Ss \subset \A$ generates the $\sigma$-field $\A$ and is closed under the formation of finite intersections (i.e.\ is a $\pi$-system), then $\Ss$ is a separating class for $\A$, cf.\ e.g.\ \cite{Bill99}, page 9.

\bl \label{sepcla2}
If $\Ss_1$ and $\Ss_2$ are separating classes for the $\sigma$-fields $\A_1$ and $\A_2$, respectively, then so is $\Ss_1 \times \Ss_2$ for $\A_1 \otimes \A_2$.
\el
\bop
We have to show that the two properties, $\pi$-system and generating class, extend from the marginals to the product. The former is apparent, for the latter see e.g.\ \cite{Bauer1992}, Theorem 22.1, page 151.
\eop

For any $T_0 \subset \R$ let
\[
    \FT = \big\{\pi^{-1}_T(A) \big|\, A \in \B(\R^{|T|}), T \subset T_0, |T| < \infty  \big\}
\]
and
\[
    \mathrsfs{H}(T_0) = \big\{ \piTT^{-1}(A) \big|\, A \in \B(\R^{2 |T|})\ , T \subset T_0, |T| < \infty 
     \big\}.
\]
$\FT$ and $\mathrsfs{H}(T_0)$ are subclasses of $\D$ and $\D \otimes \D$, respectively, cf.\ \cite{Bill99}, Theorem 16.6. Roughly, the next two lemmas tell that $\FT$ and $\mathrsfs{H}(T_0)$ are ``large enough'', if $T_0$ is ``large enough''.

\bl \label{sepcla3}
If $T_0$ is dense in $\R$, then $\FT$ is a separating class for $\D$.
\el
\bop
See \cite{Bill99}, page 170, Theorem 16.6.
\eop

\bl \label{Hsepcla}
If $T_0$ is dense in $\R$, then $\mathrsfs{H}(T_0)$ is a separating class for $\D \otimes \D$.
\el
\bop
By lemmas \ref{sepcla2} and \ref{sepcla3}: $\FT \times \FT$ is a separating class for $\D \otimes \D$. It remains to see: $\FT \times \FT \subset \mathrsfs{H}(T_0)$. Towards this end we introduce the class
\[
    \G(T_0) = \big\{ \piST^{-1}(A) \big|\, A \in \B(\R^{|S|+|T|}), \ S, T \subset T_0, \ |S|, |T| < \infty 
     \big\}.
\]
Evidently $\FT \times \FT \subset \G(T_0)$. Furthermore it holds $\G(T_0) = \mathrsfs{H}(T_0)$. This is because any set $\piST^{-1}(A) \in \mathrsfs{H}(T_0)$ can also be written as $\pi_{T \cup S, T \cup S}^{-1}(C)$ for an apropriate set $C \subset \R^{2 |T \cup S|}$.
\eop

This concludes the preliminaries, and we present the 

\medskip
\bopt{ of Theorem \ref{wc2}}
We have, $\{ \p{1}_n \}$ and $\{ \p{2}_n \}$ are both tight, hence $\{ P_n = (\p{1}_n, \p{2}_n) \}$ is tight (Lemma \ref{Z tight}). By Prokhorov's theorem, $\{ P_n \}$ is relatively compact. (Here it is important that $\D \otimes \D$ coincides with the Borel-$\sigma$-field on $D \times D$, cf. Lemma \ref{sigma-fields}.) To each subsequence $\{ P_{n'} \}$ exists a further (sub-)subsequence $\{ P_{n''} \}$ that converges, i.e.\ there is a probability measure $Q = (Q^{(1)}, Q^{(2)})$ on $(D \times D, \D \otimes \D)$ (which of course depends on $\{ n'' \}$), such that
\[
    P_{n''} \cid Q.
\]
Lemmma \ref{Paecontinuous} allows us to apply the CMT:
\[
    P_{n''} \circ \piTT^{-1} \cid Q \circ \piTT^{-1} \qquad \mbox{ for all finite } T \subset T_{Q^{(1)}} \cap T_{Q^{(2)}}.
\]
On the other hand, \ref{wc2} (\ref{wc2 b2}) implies
\[
    P_{n''} \circ \piTT^{-1} \cid P \circ \piTT^{-1} \qquad \mbox{ for all finite } T \subset T \subset T_{P^{(1)}} \cap T_{P^{(2)}}.
\]
This means, if we let $T_0 = T_{Q^{(1)}} \cap T_{Q^{(2)}} \cap T_{P^{(1)}} \cap T_{P^{(2)}}$, then $P$ and $Q$ agree on the class $\mathrsfs{H}(T_0)$. The set $T_0$ is dense in $\R$ (corollary of Lemma \ref{camc}), thus $\mathrsfs{H}(T_0)$ is a separating class for $\D \otimes \D$ (Lemma \ref{Hsepcla}), hence $P = Q$.

Thus we know, all subsequences of $\{ P_n \}$ contain a weakly convergent sub-subsequence, and all of these sub-subsequences converge to the same limit $P$. It follows: $P_n$ converges weakly to $P$, cf.\ \cite{Bill99}, Theorem 2.6, page 20.
\eop

\section{Proof of Theorem \ref{main result}} \label{proof}

\bopt{ of Theorem \ref{main result}}
By applying Theorem \ref{wc2} the proof comes down to showing 
\begin{enumerate}[(A)]
\item \label{A} $\{ \AL \}$ is tight, 
\item \label{B} $\{ \BA \}$ is tight, and
\item \label{C} the fidis of $(\AL,\BA)$ converge to those of $(B_1^F,\Nrho)$ in the sense of \ref{wc2} (\ref{wc2 b2}). 
\end{enumerate}

Part (B) is an immediate corollary of point-process tightness criterion \ref{tc3}. The result (A) does not follow equally straightforward: In a first step we use the moment-type criterion \ref{tc2} to show that it holds for continuous $F$ (Lemma \ref{al cont tight}), and then, building on that, show it for arbitrary distribution functions $F$ (Proposition \ref{al tight}). Part (\ref{C}) is subject of Proposition \ref{al,ba fidis}. 
\eop

\bl \label{al cont tight}
If $F$ is continuous, the series of random variables $\{ \AL \}$ in $\DD$ is tight.
\el
\bop
Some straightforward calculations yield
\be \label{14.9}
    \E \Big( | \alpha_n(s) - \alpha_n(r) |^2 | \alpha_n(t) - \alpha_n(s) |^2 
         \Big) \le 6 \big( F(t)-F(r)\big)^2
\ee
for all $ r < s < t$ and all $n \in \N$, cf.\ \cite{Bill99}, page 150. Now we apply Proposition \ref{tc2}. Condition \ref{tc2} (\ref{tc2b1}) is fullfilled as a corrollary of Proposition \ref{al,ba fidis} or as a simple exercise using the multivariate CLT. Condition \ref{tc2} (\ref{tc2b2}) follows from (\ref{14.9}) with $a = b = 2$ and $H = \sqrt{6}F$.\\
\eop

\bl \label{qt a e cont}
Let $F$ be an arbitrary cdf. The quantile transformation $Q_F : D[0,1] \to D(\R) : x \mapsto x \circ F$ (i.e.\ $Q_F(x)$ is the function $t \mapsto x(F(t))$) is $\mathrsfs{L}(B_0)$-a.e.\ continuous, where $\Law(B_0)$ denotes the distribution of the Brownian bridge $B_0$.
\el
\bop
We show: $Q_F$ is continuous as a function from $(C[0,1],||\cdot||)$ to $D(\R)$. This suffices because $\Pb(B_0 \in C[0,1]) = 1$ and the Skorokhod topology coincides with the uniform topology (the one induced by the sup-norm $||\cdot||$) on $C\I$, cf.\ \cite{Bill99}, page 124. 
Now let $x_n \to x$ in $(C\I,||\cdot||)$. Since the image of $F$ is a subset of $\I$, this implies
\[
    \sup_{t \in \R} | \, x_n(F(t)) - x(F(t)) | \to 0
\]
Hence (\ref{cc1}) is satisfied with $\lambda_n \equiv \id$ (the identy function on $\R$) and we have $Q_F(x_n) \to Q_F(x)$ in $D(\R)$. 
\eop

\bp \label{al tight}
The series of random variables $\{ \AL \}$ in $\DD$ is tight for all distribution functions $F$.
\ep

\bop
The cdf of the uniform$(0,1)$ distribution, which we want to call $G$, is continuous. Thus by \ref{al cont tight}, \ref{al,ba fidis} and \ref{wc1} we have
\[
    \al_n^G \cid B_1^G \quad \mbox{in} \quad \DD.
\]
If we restrict these processes to the time domain $\I$, the convergence remains true, cf.\ \cite{Bill99}, page 174, Theorem 16.7. Furthermore, for any cdf $F$, $\AL$ and $B_1^F$ are their respective quantile transformations (cf. \ref{qt a e cont}) w.r.t.\ $F$. Hence the previous lemma allows us to apply the CMT:
\[
    \AL \cid B_1^F \quad \mbox{in} \quad \DD.
\]
A convergent sequence is tight.
\eop


\bp \label{al,ba fidis}
Let $\tau \in \R$ and $F$ be an arbitrary cdf such that Condition C.1 is satisfied. Then
\be \nonumber 
            \big( \AL(t_1), ..., \AL(t_k), \BA(t_1), ...,\BA(t_k) \big)
             \cid           
            \big( B_1^F(t_1), ..., B_1^F(t_k), \Nrho(t_1), ...,\Nrho(t_k) \big)
\ee
holds true for all $k \in \N$ and $t_1, ...,t_k \in \R$.
\ep

A few remarks before we come to the proof: Showing $\al_n \fdc B_1$ is a straightforward application of the multivariate CLT. The result $\ba_n \fdc N_0$ is also easy to get using a Poisson-type limit theorem. \it Of course\rm, it does not suffice to show these two statements separately. The two sequences are not independent of each other, and we need to show the finite-dimensional convergence of the joint sequence. Obviously neither of the approaches for the marginals works here.
We prove Proposition \ref{al,ba fidis} by showing the pointwise convergence of the corresponding characteristic functions. This includes some lengthy calculations, so we restrict our demonstration to 
\[
    \big(\al_n(t), \ba_n(t)\big) \cid \big(B_1(t),N_0(t)\big),
\]
and only for $t \ge 0$. The case $t < 0$ works just the same. The full-detail proof treating arbitrary tuples $t_1, ..., t_k \in \R$ is written down in \cite{Vogel2005}. For the characteristic function of $(\al_n(t), \ba_n(t))$ we write $\psi_n^{(t)}$, or short $\psi_n$, and $\psi^{(t)}$ or $\psi$ for the \cf\ of $(B_1(t),N_0(t))$. Since $B_1(t)$ and $N_0(t)$ are independent, we can write down, if $t \ge 0$,
\be \label{psi}
    \psit(x,y) = \exp \big\{ -\frac{1}{2} F(t)(1-F(t)) x^2 + \varrho_2 t (e^{i y} -1) \big\}.
\ee
If $t < 0$, then $\psit(x,y)$ contains a term with $\varrho_1$ instead of $\varrho_2$. The next lemma specifies $\psit_n$.

\bigskip
\bl \label{psin}
Assume $t \ge 0$.
\begin{enumerate}[(1)]
\item \label{psin1}
    If $t > \tau$, and $n$ is sufficiently large such that $\tau + \frac{t}{n} < t$, then the \cf\ $\psint:\R^2 \to \C$ \ of \ $\big(\AL(t),\BA(t)\big)$ is given by \\
    $\psint(x,y) = \exp\{ -i x \sqrt{n} F(t) \} 
        \bigg[ 
            1 + F(t)\big( \gross{e^{\frac{i x}{\sqrt{n}}}}-1 \big) 
            + \big(\Ftautn -F(\tau) \big) (e^{i y} -1)\gross{e^{\frac{i x}{\sqrt{n}}}}
        \bigg]^n$.
\item \label{psin2}
    If $t \le \tau$, then $\psint:\R^2 \to \C$ is \\
    $\psint(x,y) = \exp\{ -i x \sqrt{n} F(t) \} 
        \bigg[ 
            1 + F(t)\big( \gross{e^{\frac{i x}{\sqrt{n}}}}-1 \big) 
            + \big(\Ftautn -F(\tau) \big) (e^{i y} -1)
        \bigg]^n$.
\end{enumerate}
\el
\bop 
Keep in mind that $t \ge0$. Per definition of the \cf,
\[
    \psi_n(x,y) = \E \exp \Bigg\{ 
        i x \bigg[ \frac{1}{\sqrt{n}} \sum\limits_{k=1}^n \big( \Varind{(-\infty,t]}(X_k) - F(t)     \big) \bigg] 
    +   i y \bigg[                       \sum\limits_{k=1}^n        \Varind{(\tau,\tau+\frac{t}{n}]}(X_k)           \bigg]  
    \Bigg\}.
\]
Since the $X_k$, $k = 1,...,n$, are i.i.d., this transforms to 
\be \label{psin3}
    \big( \psi_n(x,y) \big)^{\frac{1}{n}} 
    = \E \exp \bigg\{ 
        \frac{i x}{\sqrt{n}} \big( \varind{(-\infty,t]}(X_1) - F(t)     \big)
        + i y \varind{(\tau,\tau+\frac{t}{n}]}(X_1)
    \bigg\}.
\ee
The right-hand side is the expectation over a function of the discrete random variable 
$\big( \, \varind{(-\infty,t]}(X_1), \, \varind{(\tau,\tau+\frac{t}{n}]}(X_1) \, \big)$, the distribution of which we know.
\[
{\renewcommand{\arraystretch}{1.5}
\begin{array}{c@{\ }|c@{\ }|c@{\ }}   
\mbox{value of} & \multicolumn{2}{c}{\mbox{corresponding probability if}} \\
\big( \varind{(-\infty,t]}(X_1), \varind{(\tau,\tau+\frac{t}{n}]}(X_1) \big) & \tau+\frac{t}{n} < t & t \le \tau \\
\hline
(\ 0\ ,\ 0\ ) & 1 - F(t)                        & 1 - \Ftautn + F(\tau) - F(t)  \\
(\ 0\ ,\ 1\ ) & 0                               & \Ftautn - F(\tau)                     \\
(\ 1\ ,\ 0\ ) & F(t) - \Ftautn + F(\tau)    & F(t)                                  \\
(\ 1\ ,\ 1\ ) & \Ftautn - F(\tau)           & 0                                     \\
\end{array}
} 
\]

Thus in both cases we can write down the expectation (\ref{psin3}) as a sum of three summands. Some re-grouping yields the expressions in Lemma \ref{psin}.
\eop

\bopt{ of Proposition \ref{al,ba fidis}}
We show $\psint \to \psit$ pointwise  for all $t \ge 0$. We apply the following result from complex analysis. For complex numbers $c$ and $c_n$, $n \in \N$,
\be \label{ca}
    c_n \lra c \quad \Longrightarrow \quad \Big(1 + \frac{c_n}{n}\Big)^n \lra e^c.
\ee
Hence, it suffices to prove
\be \label{psin4}
    n \Big( \psi_n^{\frac{1}{n}} (x,y) -1 \Big) \lra \ln \psi(x,y).
\ee
Call the left-hand side $h_n$ and the right-hand side $h$. Consider at first the case \ref{psin} (\ref{psin1}), i.e.\ $\tau < t$. Then by (\ref{psi}) and Lemma \ref{psin},
\begin{eqnarray*}
h_n & = & n \exp \Big\{ - \frac{i x }{\sqrt{n}} F(t) \Big\} 
        \bigg[ 
            1 + F(t)\big( \gross{e^{\frac{i x}{\sqrt{n}}}}-1 \big) 
            + \big(\Ftautn -F(\tau) \big) (e^{i y} -1)\gross{e^{\frac{i x}{\sqrt{n}}}}
        \bigg] - n, \\
h & = &  -\frac{1}{2} F(t)(1-F(t)) x^2 + \varrho_2 t (e^{i y} -1).
\end{eqnarray*}
We break the convergence $h_n \to h$ down into two parts:
\begin{enumerate}[(a)]
\item $n \exp\Big\{ - \frac{i x }{\sqrt{n}} F(t) \Big\}     
        \bigg[ 1 + F(t)\big( \gross{e^{\frac{i x}{\sqrt{n}}}}-1 \big)
        \bigg] - n
        \ \lra \ 
        -\frac{1}{2} F(t)(1-F(t)) x^2$,
\item $n \exp \Big\{ \frac{i x }{\sqrt{n}} \big(1-F(t)\big) \Big\}  
        \Big[ \big(\Ftautn -F(\tau) \Big] (e^{i y} -1)
        \ \lra \
        \varrho_2 t (e^{i y} -1)$.
\end{enumerate}
For (a): Use the Taylor expansion of the exponential function. Bear in mind that it converges uniformly on any compact set. This allows us to write
\[
\exp \Big\{- \frac{i x F(t) }{\sqrt{n}} \Big\} = 1 - \frac{i x F(t)}{\sqrt{n}}  - \frac{x^2 F(t)^2}{2 n} + o\big(\frac{1}{n}\big) \qquad (n \to \infty)
\]          
and 
\[
\exp \Big\{ \frac{i x }{\sqrt{n}} \Big\} = 1 + \frac{i x }{\sqrt{n}} + o\big(\frac{1}{\sqrt{n}}\big) \qquad (n \to \infty).
\]  
Plug this into the left-hand side, the rest is computing. 

Part (b) becomes apparent by noting 
\[
    \lim\limits_n \frac{n}{t} \Big( \Ftautn -F(\tau) \Big) = \varrho_2  \qquad (t \ge 0)
\]
and
\[
    \lim\limits_n \exp \Big\{ \frac{i x }{\sqrt{n}} \big(1-F(t)\big) \Big\} = 1.
\]

By adding (a) and (b) we have proved $h_n \to h$. By (\ref{ca}) this implies $\psint(x,y) \to \psit(x,y)$, but so far only for $\tau < t$. As for $\tau \ge t$, the only difference is that in (b) the term $\ \exp\! \big\{ \frac{i x }{\sqrt{n}} (1-F(t)) \big\}$\, is replaced by $\ \exp\! \big\{- \frac{i x }{\sqrt{n}}F(t) \big\}$, \,which of course converges to $1$ as well.
\eop

\section{Application in statistics}

In this short section we demonstrate at an example how Theorem \ref{main result} can be useful in statistics. Our arguments will only briefly be sketched. Consider i.i.d.\ random variables $X_1,...,X_n$, $n \in \N$, with values in $\I$ and common cdf $F = F_{\tau,\gamma}$. Here, $F_{\tau,\gamma}$ is defined on $\I$ as the polygonal line through the points $(0,0)$, $(\tau,\gamma)$ and $(1,1)$, where the parameters $\tau$ and $\gamma$ both lie in the open interval $(0,1)$, and it is assumed that $\tau \neq \gamma$. Thus, $\tau$ is the single point of discontinuity of the corresponding density. In this model Chernoff and Rubin \cite{Chernoff1956} investigate the maximum likelihood estimator for $\tau$. An ad hoc estimator for the two-dimensional parameter $(\tau,\gamma)$ is given by
\[
	\hat{\tau}_n = \argmax_{t \in \R} |\, \dsF_n(t) - t | \quad \mbox{ and } \quad \hat{\gamma_n} = \dsF_n(\hat{\tau}_n).
\]
A key role in the analysis of the pair 
\[
	\Big( n(\hat{\tau}_n - \tau), \ \sqrt{n} (\hat{\gamma}_n - \gamma) \Big)
\]
plays the observation that it has the same limit distribution as
\[
	\Big(\, \argmax_{t \in \R} \big\{ \sign(\gamma-\tau)\big(\BA(t)-t\big) \big\},\ \AL(\tau) \, \Big).
\]
Thus Theorem \ref{main result} and a (formal) application of the CMT yield convergence in distribution:
\be \label{applCID}
	\Big( n(\hat{\tau}_n - \tau), \ \sqrt{n} (\hat{\gamma}_n - \gamma) \Big) \cid (A,B),
\ee
where
\[
	A = \argmax_{t \in \R} \big\{ \sign(\gamma-\tau)\big(N_0(t)-t\big)\big\}
\]	
and $B \sim N\big(0, F(\tau)(1- F(\tau))\big)$ are independent. The rates of $N_0 = N_0^{\varrho_1,\varrho_2}$ are given by $\varrho_1 = \gamma/\tau$ and $\varrho_2 = (1-\gamma)/(1-\tau)$. In \cite{Ferger2005} we give a representation of $A$ in terms of arrival times of $N_0$, which shows that with probability 1 the maximizing point $A$ is uniquely determined. Moreover $A$ is seen to have a continuous cdf. A rigorous proof of (\ref{applCID}) and further information will be published elsewhere.


\begin{thebibliography}{vdVW96}

\bibitem[AHE84]{Al-Hussaini1984}
A.~Al-Hussaini and R.~J. Elliot.
\newblock Convergence of the empirical distribution to the {P}oisson process.
\newblock {\em Stochastics}, 13:299--308, 1984.

\bibitem[Bau92]{Bauer1992}
H.~Bauer.
\newblock {\em Ma{\ss}- und Integrationstheorie}.
\newblock De Gruyter Lehrbuch. Walter de Gruyter \& Co., 2nd edition, 1992.

\bibitem[Bil99]{Bill99}
P.~Billingsley.
\newblock {\em Convergence of Probability Measures}.
\newblock Wiley Series in Probability and Statisitics. John Wiley \& Sons, 2nd
  edition, 1999.

\bibitem[CH88]{Csorgo1988}
M.~Cs{\"o}rg{\H{o}} and L.~Horv{\'a}th.
\newblock {Convergence of the empirical and quantile distributions to Poisson
  measures.}
\newblock {\em Stat. Decis.}, 6:129--136, 1988.

\bibitem[CR56]{Chernoff1956}
H.~Chernoff and H.~Rubin.
\newblock {The estimation of the location of a discontinuity in density.}
\newblock {\em {Proc. 3rd Berkeley Sympos. Math. Statist. Probability}},
  1:19--37, 1956.

\bibitem[Don51]{Donsker1951}
M.~D. Donsker.
\newblock An invariance principle for certain probability limit theorems.
\newblock {\em Mem. Amer. Math. Soc.}, 6, 1951.

\bibitem[Don52]{Donsker1952}
M.~D. Donsker.
\newblock Justification and extension of {D}oob's heuristic approach to the
  {K}olmogorov-{S}mirnov theorems.
\newblock {\em Ann. Math. Stat.}, 23:277--281, 1952.

\bibitem[Doo49]{Doob1949}
J.~L. Doob.
\newblock Heuristic approach to the {K}olmogorov-{S}mirnov theorems.
\newblock {\em Ann. Math. Stat.}, 20:393--403, 1949.

\bibitem[EK86]{Ethier1986}
S.~N. Ethier and T.~G. Kurtz.
\newblock {\em Markov Processes: Characterization and Convergence}.
\newblock John Wiley \& Sons, 1986.

\bibitem[Els02]{Elstrodt2002}
J.~Elstrodt.
\newblock {\em Ma{\ss}- und Integrationstheorie}.
\newblock Grundwissen Mathematik. Springer, 3rd edition, 2002.

\bibitem[Fer05]{Ferger2005}
D.~Ferger.
\newblock {On the minimizing point of the incorrectly centered empirical
  process and its limit distribution in nonregular experiments.}
\newblock {\em ESAIM, Probab. Stat.}, 9:307--322, 2005.

\bibitem[Gal78]{Galambos1978}
J.~Galambos.
\newblock {\em {The asymptotic theory of extreme order statistics.}}
\newblock Wiley Series in Probability and Statisitics. {John Wiley \& Sons},
  1978.

\bibitem[JS02]{Jacod2002}
J.~Jacod and A.~N. Shiryaev.
\newblock {\em Limit Theorems for Stochastic Processes}.
\newblock Springer, 2nd edition, 2002.

\bibitem[KLS80]{Kabanov1980}
Yu.~M. Kabanov, R.~Sh. Liptser, and A.~N. Shiryaev.
\newblock Some limit theorems for simple point processes (a martingale
  approach).
\newblock {\em Stochastics}, 3:203--216, 1980.

\bibitem[Kol33]{Kolmogorov1933}
A.~N. Kolmogorov.
\newblock {Sulla determinazione empirica di una legge di distribuzione.}
\newblock {\em Giorn. Ist. Ital. Attuari}, 4:83--91, 1933.

\bibitem[Kol56]{Kolmogorov1956}
A.~N. Kolmogorov.
\newblock On {S}korokhod convergence.
\newblock {\em Theory Probab. Appl.}, 1:215--222, 1956.

\bibitem[Lin73]{Lindvall1973}
T.~Lindvall.
\newblock Weak convergence of probability measures and random functions in the
  function space ${D}[0,\infty)$.
\newblock {\em J. Appl. Probab.}, 10:109--121, 1973.

\bibitem[Pol84]{Pollard1984}
D.~Pollard.
\newblock {\em Convergence of Stochastic Processes}.
\newblock Springer, 1984.

\bibitem[Pro56]{Prokhorov1956}
Yu.~V. Prokhorov.
\newblock Convergence of random processes and limit theorems in probability
  theory.
\newblock {\em Theory Probab. Appl.}, 1:157--214, 1956.

\bibitem[Ros67]{Rossberg1967}
H.~J. Rossberg.
\newblock {\"Uber das asymptotische Verhalten der Rand- und Zentralglieder
  einer Variationsreihe. II.}
\newblock {\em Publ. Math.}, 14:83--90, 1967.

\bibitem[RY99]{Revuz1999}
D.~Revuz and M.~Yor.
\newblock {\em {Continuous martingales and Brownian motion}}.
\newblock {Grundlehren der Mathematischen Wissenschaften}. {Springer}, 3rd
  edition, 1999.

\bibitem[Sko56]{Skorokhod1956}
A.~V. Skorokhod.
\newblock Limit theorems for stochastic processes.
\newblock {\em Theory Probab. Appl.}, 1:261--290, 1956.

\bibitem[Smi44]{Smirnov1944}
N.~V. Smirnov.
\newblock {Approximate laws of distribution of random variables from empirical
  data (in Russian).}
\newblock {\em Usp. Mat. Nauk}, 10:179--206, 1944.

\bibitem[Sto63]{Stone1963}
C.~Stone.
\newblock Weak convergence of stochastic processes defined on a semi-infinite
  time interval.
\newblock {\em Proc. Amer. Math. Soc.}, 14:694--696, 1963.

\bibitem[SW86]{Shorack1986}
G.~R. Shorack and J.~A. Wellner.
\newblock {\em Empirical Processes with Applications to Statistics}.
\newblock Wiley Series in Probability and Mathematical Statistics. John Wiley
  {\&} Sons, 1986.

\bibitem[vdVW96]{Vaart1996}
A.~W. van~der Vaart and J.~A. Wellner.
\newblock {\em Weak Convergence and Empirical Processes}.
\newblock Springer Series in Statistics. Springer, 1996.

\bibitem[Vog05]{Vogel2005}
D.~Vogel.
\newblock Weak convergence of the empirical process and the rescaled empirical
  distribution function in the {S}korokhod produkt space.
\newblock Master's thesis, Technische Universit{\"a}t Dresden, Dresden,
  Germany, November 2005.

\bibitem[Whi70]{Whitt1970}
W.~Whitt.
\newblock Weak convergence of probability measures on the function space
  {C[0,1]}.
\newblock {\em Ann. Math. Stat.}, 41:939--944, 1970.

\bibitem[Whi71]{Whitt1971}
W.~Whitt.
\newblock Weak convergence of probability measures on the function space
  {$D[0,\infty)$}.
\newblock Technical report, Yale University, 1971.

\bibitem[Whi02]{Whitt2002}
W.~Whitt.
\newblock {\em Stochastic-Process Limits}.
\newblock Springer Series in Operations Research. Springer, 2002.

\end{thebibliography}
\end{document}